\documentclass[a4paper,10pt]{article}

\usepackage{bm}
\usepackage{bbm}
\usepackage[utf8]{inputenc}
\usepackage[english]{babel}
\usepackage{amssymb}
\usepackage{amsmath}
\usepackage{amsthm}
\usepackage{amsfonts}
\usepackage{mathrsfs}
\usepackage{dsfont}
\usepackage{hyperref}
\usepackage{commath}
\usepackage{graphicx}
\usepackage{text comp}
\usepackage{authblk}
\usepackage{mathtools}
\usepackage{afterpage}
\usepackage{appendix}

\theoremstyle{plain}
\newtheorem*{thm*}{Theorem}
\newtheorem{thm}{Theorem}[section]

\newtheorem{lemma}[thm]{Lemma}
\newtheorem*{lemma*}{Lemma}

\newtheorem*{corollary*}{Corollary}

\newtheorem*{prop*}{Proposition}
\newtheorem{defn}{Definition}

\newtheorem*{conjecture*}{Conjecture}

\newcommand{\R}{\mathbb{R}}

\newcommand{\N}{\mathbb{N}}

\newcommand{\1}{\mathbbm{1}}

\newcommand{\dist}{\mathrm{dist}}

\everymath{\displaystyle}

\title{Random \v Cech complexes on $\R^d$: decrackling the noise with local scalings}
\author[ ]{Henry-Louis de Kergorlay}
\affil[ ]{School of Mathematics, University of Edinburgh}

\date{}

\begin{document}

\maketitle
\begin{abstract}
    We investigate the homology of an unbounded noisy sample on $\R^d$, under various assumptions on the sampling density. This investigation is based on previous results by Adler, Bobrowski, and Weinberger (\cite{crackle}), and Owada and Adler (\cite{topoCrackle}). There, it was found that unbounded noise generally introduces non-vanishing homology, a phenomenon called \textit{topological crackle}, unless the density has superexponential decay on $\R^d$. We show how some well-chosen \textit{non-trivial} variable bandwidth constructions can extend the class of densities where crackle doesn't occur to any light tail density with mild assumptions, what we call \textit{decrackling the noise}.
\end{abstract}
\section{Introduction}
A random geometric complex is a simplicial complex built on a random geometric graph. Topological properties of random geometric complexes have been investigated (e.g., \cite{randGeom,vanishHomo,randCech, survey,cechBoundary}) as generalizations of properties of random geometric graphs, as studied by Penrose in \cite{randGraphs}.\\

As observed in \cite{randGeom}, the investigation of topological properties of random geometric complexes is motivated by applications in Topological Data Analysis. Investigating topological features of randomly sampled points can serve as a null hypothesis for statistical persistent homology (e.g., \cite{computPersHom,persHomSurvey,TDA,statPersistent,persistentCycles}). In particular, it is desirable to find asymptotic bounds on the expected Betti numbers (the expected ranks of the homology groups) under various settings, and to exhibit\  connectivity thresholds (asymptotic values of the bandwidth beyond which homology vanishes w.h.p.).\\

The two most commonly studied complexes are the \v Cech complex and the Vietoris-Rips complex. The \v Cech complex has the advantage of yielding a geometric interpretation via the Nerve Lemma, making it the more natural complex to study when points are sampled from a Riemaniann submanifold $M\subset\R^d$ (e.g., \cite{vanishHomo,randCech,cechBoundary}). The Vietoris-Rips complex on the other hand, can be regarded as a completely combinatorial complex - it is the clique complex of a random graph (also known as the flag complex).\\

Expected topological features of the \v{C}ech and the Vietoris-Rips complex can be investigated for various phases of the bandwidth parameter $r$.
As long as the bandwidth satisfies $r= O(n^{-1/d})$ where $n$ is the number of vertices, i.e., that we are in the subcritical or critical phase, topological investigations of complexes can be done with fairly general assumptions on the sampling density, which choice does not affect greatly the asymptotic expected values of the Betti numbers (see \cite{randGeom,survey}). On the other hand, as observed in \cite{randGraphs,randGeom,survey}, a change in the sampling density may significantly alter the resulting connectivity threshold of a graph. In particular in the supercritical phase, i.e., $r=\omega(n^{-1/d})$, one generally restricts the density for instance to being compactly supported and bounded away from $0$, in order to identify a connectivity threshold value for $r$ (beyond which the homology groups vanish w.h.p.).\\


It is shown in \cite{randGeom}  that for a uniform sampling density on a smooth convex body of $\R^d$, the expected Betti numbers of the \v Cech and the Vietoris-Rips complex grow sublinearly in the supercritical phase. It is also shown that the homology of the \v Cech complex vanishes beyond the threshold value of $r\sim (\log n)^{1/d}n^{-1/d}$. In order to identify such a threshold value for the \v Cech complex, it suffices by the Nerve Lemma to identify a value beyond which the union of the balls $\cup_{x\in X_n}B(x,r)$ is contractible.\\ 

Such results are classic and well known by now; they have been expanded in various ways. The following list of works is by no means exhaustive.\\ Bobrowski and Weinberger \cite{vanishHomo} and later Bobrowski and Oliveira \cite{randCech} continued the investigation of \v Cech homology with more generality assumed on the sampling domain. They consider respectively the case where the domain is a torus and more generally when it is a compact and closed (without boundary) Riemannian manifold. The task in those settings is to identify upper and lower thresholds such that if the bandwidth parameter $r$ is below the lower threshold then w.h.p. the homology of the \v Cech complex does not recover that of the underlying manifold, while if $r$ is beyond the upper threshold it does so w.h.p..\ The upper and lower threshold have tight gaps and these sharp transitions of states are reminiscent of the well known sharp transition threshold value for the connectivity of a random geometric graph, given by $\Lambda:=n\omega_d r^d=\log n$, where $\omega_d$ denotes the Euclidean volume of the unit ball in $\R^d$ ($d$ being the dimension of the manifold). See also the recent work of Bobrowski in \cite{homologicalConnectivity} for sharp homological connectivity threshold values for the random \v Cech complex, occuring at $\Lambda=\log n+(k-1)\log\log n$. The case of a compact Riemannian manifold with non-empty boundary was recently investigated in \cite{cechBoundary}. Upper and lower threshold values were provided as in the case of a closed compact manifold, but with a different scaling due to the presence of a non-empty boundary, complicating several steps. Although the difference between the two thresholds remains small, no sharp transition threshold has yet been identified in this setting. In all of the above mentioned works, only the case of a uniform density was considered. Note (as observed in \cite{survey}) that by compactness, these results easily generalize to an arbitrary density with the extra requirements that it remains bounded and strictly positive on the domain.\\

In this paper we focus our topological investigation on the \v Cech complex where the sampling density has support $\R^d$. This setting has already been carefully investigated in \cite{crackle} and \cite{topoCrackle}. In \cite{crackle} the authors are interested in the effect of unbounded noise on homology recovery of a bounded domain. They take the bounded domain to be just the origin and model the unbounded noise by sampling densities supported on $\R^d$, with various decays (\textit{power law}, \textit{exponential}, \textit{Gaussian}). Thus, they identify sequences of radii $(R_n^c)_{n\in \N}$ depending on the density, such that homology vanishes in $B(0,R_n^c)$ but not necessarily outside the ball, in which case they call this phenomenon \textit{topological crackle}. If topological crackle occurs for a given density, then cycles still form far away from the origin, even for $r_n>> (\log n)^{1/d}n^{-1/d}$, i.e., the homology of the noise doesn't vanish and we cannot hope to recover the true homology from a noisy sample. In \cite{crackle} the authors show that crackling occurs for the power law and the exponential decay, but not for the Gaussian. In \cite{topoCrackle} the authors pursue the work initiated in \cite{crackle} with more generality assumed on the density, using tools from Extreme Value Theory. They exhibit finite sequences of annuli splitting $\R^d$ (as in \cite{crackle}), where each annulus generates homology of a certain degree. Their analysis confirms that topological crackle generally occurs (in which case, as mentioned above, we cannot recover the true homology from a noisy sample) unless the density has superexponential decay, a generalization of the Gaussian case analyzed in \cite{crackle} (see Theorem $4.5$ in \cite{topoCrackle}).\\

Our main results, summarized in Theorem \ref{main results}, provide well-chosen local scalings and asymptotic conditions on the bandwidth parameter $r_n$ depending on the sampling density, such that the homology groups of the associated variable bandwidth \v Cech complexes (see Definition \ref{variable bandwidth Cech}) vanish with high probability (w.h.p.). The novelty of these results compared to those in \cite{crackle,topoCrackle}, is that they provide well-chosen variable bandwidth constructions which allow us to identify asymptotic conditions on the bandwidth parameter in order for the homology of the noise (of the \v Cech complex) to vanish in cases previously shown to induce topological crackle. Indeed, while suitable conditions on the bandwdith parameter in order to avoid the crackling phenomenon can only be found in the case where the sampling density has Gaussian or more generally superexponential decay following the results in \cite{crackle,topoCrackle}, our constructions will allow us to identify such conditions on the bandwidth parameter for any light tail distribution under mild assumptions. Furthermore, if the density has  superexponential decay we present a construction which can weaken the asymtpotic conditions on the bandwith parameter. \\


\section*{Acknowledgements}
I would like to acknowledge support from the EPSRC Programme Grant\\ EP/P020720/1. I also like to acknowledge support from the EPSRC studentship EP/L016508/1, when part of this work was initiated. I am grateful to The Alan Turing Institute for hosting me as an Enrichment student in 2017-2018, and to Ulrike Tillmann and Oliver Vipond for insightful discussions during this time, on variable bandwidth constructions and more generally on stochastic topology. I would also like to thank them and Leonardo Tolomeo for great discussions and comments on ways to extend the results below to more general unbounded manifolds. Finally, I want to thank Desmond Higham and Ulrike Tillmann who read previous versions of this work and gave many helpful comments and suggestions.

\section{Geometric complexes}

To motivate the definition of \v Cech and Vietoris-Rips complexes, recall the following definition of a geometric graph.\\

\begin{defn}[Geometric graph]
Given vertices $\mathcal P\subset \R^d$ and a radius $r>0$ (also called bandwidth parameter) define the geometric graph $G(\mathcal P;r)$ where $\{x,y\}\subset\mathcal P$ is an edge if
$$
B(x,r/2)\cap B(y,r/2)\neq\emptyset.
$$
\end{defn}

 Geometric complexes can be regarded as generalizations of geometric graphs, where we not only take into consideration vertices and edges, but also triangles and higher dimensional simplices. As mentioned above, there are two common geometric complexes, which we now define.\\
 
 \begin{defn}[\v Cech complex]
 Given vertices $\mathcal P\subset\R^d$ and bandwidth parameter $r>0$, define the \v Cech complex $\mathcal C(\mathcal P,r)$ to be the simplicial complex where for all $k\in [d]$, $[x_0,\dots,x_k]$ is a $k$-face of the complex if
$$
\cap_{i=0}^kB(x_i,r)\neq \emptyset.
$$
 \end{defn}
 \begin{defn}[Vietoris-Rips complex]
 Given vertices $\mathcal P\subset\R^d$ and $r>0$, define the Vietoris-Rips complex $R(\mathcal P,r)$ to be the clique complex of the graph $G(\mathcal P,r)$, i.e., for all $k\in \N$,  $[x_0,\dots,x_k]$ is a $k$-face of the complex if it is a clique in $G(\mathcal P,r)$, which means
 $$
 \forall i_1\neq i_2,\ B(x_{i_1},r/2)\cap B(x_{i_2},r/2)\neq \emptyset.
 $$
 \end{defn}
 
 Those two complexes are closely related. In fact, it is immediate to verify the following inclusions:
 $$
 \mathcal C(\mathcal P,r)\subset R(\mathcal P,2r)\subset \mathcal C(\mathcal P,2r).
 $$
 For a tighter nested inclusion, see \cite{cechVR}.\\
 
 While the Vietoris-Rips complex is completely combinatorial, in that all the information is contained in the underlying graph (i.e., it is the clique complex of the graph), the \v Cech complex has the advantage of yielding a nice geometric interpretation via the following celebrated Nerve Lemma. Below, we phrase it in a way directly related to our purposes. A more general and well known topological result holds, originating in \cite{nerveOrigins}.
 
 \begin{lemma}[Nerve Lemma]
 Let $\mathcal P\subset \R^d$, and let $\mathcal B(\mathcal P,r):=\{B(x,r)\ |\ x\in \mathcal P\}$. The \v Cech complex $\mathcal C(\mathcal P,r)$ is homotopy equivalent to $\mathcal B(\mathcal P,r)$.
 \end{lemma}
 
 The Nerve Lemma gives us a simple geometric criterion for the vanishing of \v Cech homology: the contractibility of the union of the balls $\cup_{x\in \mathcal P} B(x,r)$.

 \section{Local scalings: varying the bandwidth}
 
 Variable bandwidth constructions are well known and have been successfully used in various topics, for instance in non-parametric kernel density estimations (e.g., \cite{varKerEstimates,varKerEstimations,locallyAdaptive}) and spectral clustering (e.g., \cite{selfTuning}). The main idea consists in choosing a local scaling varying as a function of the vertices, as opposed to a constant scaling, in constructions involving a kernel function. Opting for a local scaling has many advantages. As discussed in \cite{selfTuning}, in the case where the sampling density is non-uniform, suitable scaling choices can allow us to handle multi-scale data.\\
 
 There are various ways one can define a so-called self-tuning or variable bandwidth geometric graph. In \cite{selfTuning} the authors choose the affinity matrix as
 $$
 \1\left(|x-y|^2\leq \sigma(x)\sigma(y)\right), 
 $$
 where $\sigma:\R^d\to \R_+^*$.
 Here we take the associated affinity matrix to be instead
 $$
\1\left(|x-y|\leq \frac{\sigma(x)+\sigma(y)}{2}\right).
$$

\begin{defn}[Variable bandwidth geometric graph]
 Given vertices $\mathcal P\subset \R^d$ and a (scaling) function $\sigma:\R^d\to \R_+^*$, define the variable bandwidth geometric graph $G(\mathcal P,\sigma)$ to be the graph with vertex set $\mathcal P$ and where $\{x,y\}\subset \mathcal P$ is an edge if
 $$ B(x,\sigma(x)/2)\cap B(y,\sigma(y)/2)\neq\emptyset.
 $$
 \end{defn}

Likewise, the definitions of the \v Cech and the Vietoris-Rips complex naturally extend to a variable bandwidth setting.
\begin{defn}[Variable bandwidth \v Cech complex]\label{variable bandwidth Cech}
Given vertices $\mathcal P\subset\R^d$ and a scaling $\sigma:\R^d\to \R_+^*$, define the \v Cech complex $\mathcal C(\mathcal P,\sigma)$ to be the simplicial complex where for all $k\in [d]$, $[x_0,\dots,x_k]$ is a $k$-face of the complex if
$$
\cap_{i=0}^kB(x_i,\sigma(x_i))\neq \emptyset.
$$
\end{defn}
\begin{defn}[Variable bandwidth Vietoris-Rips complex]
 Given vertices $\mathcal P\subset\R^d$ and a scaling $\sigma:\R^d\to \R_+^*$, define the Vietoris-Rips complex $R(\mathcal P,\sigma)$ to be the clique complex of the graph $G(\mathcal P,\sigma)$, i.e., for all $k\in \N$,  $[x_0,\dots,x_k]$ is a $k$-face of the complex if it is a clique in $G(\mathcal P,\sigma)$, which means
 $$
 \forall i_1\neq i_2,\ B(x_{i_1},\sigma(x_{i_1})/2)\cap B(x_{i_2},\sigma(x_{i_2})/2)\neq \emptyset.
 $$
\end{defn}

Picking a suitable scaling for $\sigma$ is not straightforward and different choices may yield very different results. A classic scaling (e.g., \cite{locallyAdaptive,varKerEstimations,varKerEstimates,selfTuning,variableBandwidth}) generally consists in choosing $\sigma(x_i)= \abs{x_i-x_i^{(k)}}$, where $x_i^{(k)}$ is the $k$-nearest neighbour of $x_i$ from the sampled points. Such a scaling choice has the advantage not to rely on previous knowledge of the sampling density $q$ or the ambient dimension $d$. Furhtermore as $n\to \infty$, $k\to \infty$ and $k=o(n)$, it can be be asymptotically approximated as
\begin{equation}\label{approximation}
\abs{x-x^{(k)}}\sim r_n\varphi(x),
\end{equation}
where $\varphi(x) := q(x)^{-1/d}$ and $r_n:=(k/\omega_d n)^{1/d}$, and where $\omega_d$ is the Euclidean volume of a unit ball in $\R^d$ (e.g., \cite{efficientStatVar}). One easily verifies that this local scaling makes the measure of the balls (i.e., the probability for a random point to lie in that ball) roughly constant everywhere in the domain, for fixed (and small) $r$, a property which is characteristic of the uniform distribution and which facilitates calculations in various settings.\\

A naive scaling choice for a variable bandwidth construction on $\R^d$ (or any unbounded domain) would consist in choosing similarly than in the case of a bounded domain $\sigma(x)\sim  r\varphi(x)$. If the domain is bounded, with this local scaling and the following classic assumptions on the sampling density
$$
0<q_{\min}\leq q_{\max}<\infty,
$$
then the results of Kahle in Section $6$ of \cite{randGeom}, giving contractibility threshold values for the bandwidth parameter $r$, immediately generalize to the variable bandwidth setting with the same threshold values for $r$ up to multiplicative constants, both for the \v Cech (Theorem $6.1$ in \cite{randGeom}) and the Vietoris-Rips complex (Theorem $6.5$ in \cite{randGeom}).\\

On the other hand, if the domain is unbounded this scaling choice presents major drawbacks. While it is possible to extend the argument found in \cite{randGeom} with a variable bandwidth construction when the sampling density is arbitrary and supported on all of $\R^d$, it is in fact straightforward to see that with the scaling $\sigma\sim r\varphi$ the radii of balls far from the origin must grow so fast that beyond a certain radius value all the balls contain the origin. At this point the construction becomes trivial and uninteresting, the graph being complete.\\
 
 
 Following the above discussion, we consider below cases where a more sensible scaling choice than $\sigma=r\varphi$ can be made, in order to ``decrackle" the noise in a non-trivial way (see Definition \ref{valid scaling}). Let us first introduce the setting on $\R^d$ we shall be working with.\\




\section{Setting on $\R^d$}
For $\R^d$, we relate to and follow closely the setting considered in \cite{topoCrackle}. There, sequences of radii $(R_n^c)_{n\in \N}$ and $(\overline{R_n})_{n\in \N}$ are more or less explicitly given, depending on the density, such that the union of the balls which have their centre lying in $B(0,R_n^c)$ is contractible (in fact covers the ball) w.h.p., while the probability to find points outside the ball $B(0,\overline{R_n})$ tends to $0$. Furthermore, densities with superexponential decay are shown to satisfy the extra property that $\overline{R_n}-R_n^c=o(r_n)$, under specified asymptotic conditions on the bandwidth parameter $r_n$ of the graph, such that in fact the union of the balls $\cup_{x\in \mathcal P_n}B(x,r_n)$ (not just those with centre contained in $B(0,R_n^c)$) is contractible. In this case by the Nerve Lemma, the associated \v Cech complex (that of the noise) vanishes, i.e., such that eventually there is no topological crackle. The elementary, yet key observation, is that whenever the difference $\overline{R_n}-R_n^c$ decays faster than the bandwidth parameter, one can deduce that the union of the balls is contractible and hence the \v Cech homology of the noise vanishes and there is no topological crackle (as $n\to \infty$). However, as explained above, it was shown in \cite{topoCrackle} that this desirable property happens only if the sampling density has superexponential decay. In this paper, we show how well-chosen variable bandwidth constructions allow us to extend the class of densities for which topological crackle doesn't happen. We call this \textit{decrackling the noise}. If the density has a heavy tail (cf, definition in Section $3$ in \cite{topoCrackle}), then there is no hope to find any result much more interesting than what one may obtain following the ``naive scaling approach" discussed above. This is because the radii $\overline{R_n}$ and $R_n^c$ do not have the same asymptotic order, hence $\overline{R_n}-R_n^c\sim \overline{R_n}$. This can easily be verified, for instance, in the case of the power law density
$$
q(x)=\frac{1}{1+|x|^\alpha},
$$
taking $r_n\equiv 1 $ and referring to the asymptotic values for $\overline {R_n}$ and $R_n^c$ in \cite{crackle}. In this case any variable bandwidth construction will have to be such that $r_n\sigma(R_n^c)\sim \overline{R_n}$ if we wish for the homology of the noise to vanish, so that all the balls with centre near $\partial B(0,R_n^c)$ contain the origin, which yields a trivial construction of limited interest, as discussed above.\\

In the case of light tail densities on the other hand, the key fact is that there $\overline{R_n}\sim R_n^c$, hence in particular their difference grows much more slowly than the radii themselves, i.e.,
\begin{equation}\label{precondition}
\overline{R_n}-R_n^c=o(R_n^c).
\end{equation}
This creates some room for a well-chosen scaling $\sigma$ such that
\begin{equation}\label{condition 1}
\overline{R_n}-R_n^c=o(\sigma(R_n^c)),
\end{equation}
while also
\begin{equation}\label{condition 2}
\sigma(R_n^c)=o(R_n^c).
\end{equation}

Condition (\ref{condition 1}) will ensure that the union of the balls is contractible, i.e., that the homology of the noise vanishes (w.h.p. as $n\to \infty$). Condition (\ref{condition 2}) ensures that the graph remains non-trivial (that it is not complete), i.e., that the radii of the balls far from the origin don't grow sufficiently fast to contain the origin.\\

\begin{defn}\label{valid scaling}
Given sequences $(\overline R_n)_n$ and $(R_n^c)_n$ as described above, satisfying condition (\ref{precondition}), and a (decaying) sequence of bandwidth parameters $(r_n)_n$, we say that a scaling $\sigma:\R^d\to \R_+^*$ is \textbf{noise-killing} if condition (\ref{condition 1}) is satisfied, and that it is \textbf{non-trivial} if condition (\ref{condition 2}) is satisfied. If $\sigma$ is both noise-killing and non-trivial, we say that it is \textbf{valid}.
Likewise, we call a variable bandwidth construction \textbf{valid} if the associated scaling is valid.
\end{defn}

In this paper we present some \textit{valid} variable bandwidth constructions which allow us to extend the class of densities exempt from topological crackle to any light tail density with mild extra assumptions.\\

Let us recall the definition of a light tail density given in \cite{topoCrackle}, which follows itself \cite{multiExcessDistributions, highRisks}. Let the density function be given by
\begin{align*}
q:\R^d&\longrightarrow \R_+\\
x&\longmapsto L(|x|)\exp\left(-\psi(|x|)\right),
\end{align*}
where $\psi\in C^2(\R_+;\R)$ is a function of von Mises type: such that $\forall z\to \infty,$
$$
\psi'(z)>0\text{, }\psi(z)\to \infty\text{, }a'(z)\to 0,
$$
where
$$
a(z):=\frac{1}{\psi'(z)}\text{, }z\in \R_+.
$$

Furthermore, in order to ensure that $L$ is asymptotically constant and negligible compared to the tail, we require the following conditions (see conditions $(4.2)$ and $(4.3)$ in \cite{topoCrackle}):
$$
\frac{L(t+a(t)v)}{L(t)}\to 1\text{ as } t\to \infty\text{ uniformly on intervals;}
$$
and $\exists (\gamma,z_0,C)\in \R_+\times\R_+^*\times [1,\infty),$
$$
\forall t>1,\forall z\geq z_0\text{, }\frac{L(zt)}{L(z)}\leq Ct^\gamma.
$$

In order to facilitate our calculations and have a clear statement of the results, we restrict our investigation by adding a few extra assumptions which are standard.\\
\begin{itemize}

    \item First of all, we assume that $L\equiv C$, where $C$ is a suitable normalizing constant.\\
    
    \item Secondly, we assume that $\psi\in RV_v$ for some $v>0$ (is regularly varying),
 i.e., that for every $t>0$
$$
\lim_{R\to \infty}\frac{\psi(tR)}{\psi(R)} = t^v.
$$

\end{itemize}

If $v<1$ we say that the density has subexponential decay, if $v=1$ we say that the density has exponential decay, and if $v>1$ we say that the density has superexponential decay.\\

Since $\psi$ is eventually monotone (strictly increasing), its inverse function $\psi^{\leftarrow}$ is well-defined asymptotically and is also asymptotically strictly increasing. Hence we may and do extend the domain of $\psi^\leftarrow$ to all of $\R_+$ in such a way that
$$
\forall z\in \R_+,\ 1\leq \psi^\leftarrow(z).
$$

For $n\in \N$, let $\mathcal P_n$ be a homogenous Poisson point process with intensity $n$ with respect to the density $q$, i.e.,
$$
\mathcal P_n:=\{x_1,\dots,x_N\},
$$
where the $x_i$'s are i.i.d.\ samples with respect to $q$ and independent from $N\sim Po(n)$.\\

\section{Outline of the results}

In Section $7$, we discuss the case where the light tail density has subexponential or exponential decay and show that a variable bandwidth construction can be found to decrackle the noise with a \textit{valid} scaling.\\

In Section $8$, we investigate the case of a density with superexponential decay. Even though this case is already known to be exempt from crackle, the authors in \cite{topoCrackle} asked whether the asymptotic conditions given on $r_n$ could be improved from 
\begin{equation}\label{bandwidth condition}
\frac{a\circ \psi^\leftarrow(\log n)\log \log n}{r_n}=o(1)
\end{equation}
to merely:
$$
\frac{a\circ \psi^\leftarrow(\log n)}{r_n}=o(1),
$$
where recall that $\textstyle q(x)\sim e^{-\psi(|x|)}$ and $\textstyle a(z)=\frac{1}{\psi'(z)},$ $z\in \R_+$. While this question remains open in the case of the classic ``constant" bandwidth construction investigated in \cite{topoCrackle}, we will see that a well-chosen \textit{valid} variable bandwidth construction can indeed allow us to weaken the asymptotic condition (\ref{bandwidth condition}) as asked above.\\

The main results of this paper can be summarized in the following theorem.

\begin{thm}\label{main results}
Let $q\sim \exp(-\psi)$ with $\psi\in RV_v$, for some $v> 0$.\\
\begin{enumerate}
\item If $v\leq 1$ (subexponential or exponential decay) and $(r_n)_n$ satisfies
$$
\frac{a\circ\psi^\leftarrow(\log n)\log\log n}{r_n\psi^\leftarrow(\log n)}=o(1),
$$
then the scaling
$$
x\mapsto \psi^\leftarrow(\log (e+\varphi(x)))
$$
is \textbf{valid} and
$$
\lim_{n\to \infty}\mathbb P\left(\bigcup_{x\in \mathcal P_n}B(x,2r_n\psi^\leftarrow(\log(e+\varphi(x))))\text{ is contractible }\right)=1.
$$

\item If $v>1$ (superexponential decay) and $(r_n)_n$ satisfies
$$
\frac{a\circ \psi^\leftarrow(\log n)}{r_n}=o(1),
$$
letting
\begin{align*}
    \mathcal L:\R_+&\longrightarrow \R_+\\
    z&\longmapsto \log(e+\log(1+z)),
\end{align*}
then the scaling
$$
x\mapsto \mathcal L(\varphi(x))
$$
is \textbf{valid} and
$$
\lim_{n\to \infty}\mathbb P\left(\bigcup_{x\in \mathcal P_n}B(x,2r_n\mathcal L(\varphi(x)))\text{ is contractible }\right)=1.
$$
\end{enumerate}
\end{thm}
\section{Preliminaries}

For the setting, described in Section $4$, we suppose that our light tail sampling density satisfies
$$
q(z)\sim e^{-\psi(|z|)},\ z\in \R^d,
$$
where $\psi\in RV_v$ for some $v>0$.\\

Let us start with some preliminary observations, useful for the later parts of the argument.\\

\begin{lemma}\label{preliminaries RV}
\begin{itemize}
\item It always holds that
$$
a\circ\psi^\leftarrow(z) = (\psi^\leftarrow)'(z),
$$
and that
$$
\lim_{z\to \infty}\frac{a\circ\psi^\leftarrow(z)\log (z)}{\psi^\leftarrow(z)}=0.
$$

\item If $\psi\in RV_v$ with $v\leq 1$, then
$$
\lim_{z\to \infty} a\circ \psi^\leftarrow(z)\log z =\infty.
$$
\end{itemize}
\end{lemma}

\begin{proof}
Note that
$$
1 = (\psi\circ \psi^\leftarrow)'(z) = (\psi^\leftarrow)'(z)\psi'(\psi^\leftarrow(z)),
$$
hence
$$
a\circ \psi^\leftarrow(z)=(\psi'(\psi^\leftarrow(z)))^{-1}=(\psi^\leftarrow)'(z).
$$
We have $\psi\in RV_{v}$ for some $v>0$, so $\psi^\leftarrow \in RV_{1/v}$ and 
$$
a\circ\psi^\leftarrow = (\psi^\leftarrow)'\in RV_{1/v-1}.
$$
The observations made thus far were already alluded to (without proof) in \cite{topoCrackle}. It follows that
$$
\frac{a\circ\psi^\leftarrow}{\psi^\leftarrow}\in RV_{-1},
$$
hence in particular
$$
\lim_{z\to \infty}\frac{a\circ\psi^\leftarrow(z)\log (z)}{\psi^\leftarrow(z)}=0.
$$
Finally, if $v\leq 1$ then $1/v-1\geq 0$ and as we have just seen, $a\circ\psi^\leftarrow\in RV_{1/v-1}$, hence
$$
\lim_{z\to \infty} a\circ \psi^\leftarrow(z)\log z =\infty.
$$

\end{proof}

Let us define the sequences of radii of interest $(R_n^c)_n$ and $(\overline{R_n})_n$, alluded to in the introduction, following the setup in \cite{topoCrackle}.

\begin{defn}\label{assumptions}
Let $R_n^c:=\psi^\leftarrow(A_n),$ where
$$
A_n:=\log n + d\log r_n - \log\log r_n^{-1}\psi^\leftarrow(\log n)-\delta
$$
and $\delta$ satisfies
$$
d-e^\delta g^d C <0;
$$
and let $\overline{R_n}:=\psi^\leftarrow(B_n),$ where
$$
B_n:=\log n+(d-1)\log\psi^\leftarrow(\log n)+\log a\circ\psi^\leftarrow(\log n)+\log\log n.
$$
\end{defn}
\begin{lemma}\label{critical covering}
If $\log r_n = o(\log n)$, then
$$
\lim_{n\to \infty}\mathbb P\left(B(0,R_n^c)\subset \bigcup_{\mathcal P_n\cap B(0,R_n^c)}B(x,r_n)\right)=1.
$$
\end{lemma}
This follows as in the proof of Theorem $4.5$ in \cite{topoCrackle}.
\begin{lemma}\label{crackle radius}
We have
$$
\lim_{n\to \infty}\mathbb P\left(\mathcal P_n\cap B(0,\overline{R_n})=\emptyset\right)=1.
$$
\end{lemma}
Note that this lemma is similar to part of the content of Theorem $4.5$ in \cite{topoCrackle}. We provide a proof for completeness, as we have assumed slightly weaker assumptions on $a$.

\begin{proof}
To ease notation, let us denote $R:=\overline{R_n}$ for this proof.\\
Note that
$$
\mathbb P\left(\mathcal P_n\cap B(0,R)=\emptyset\right) = \exp\left(-n\int_{|x|\geq R}q(x)dx\right).
$$
We show $n\int_{|x|\geq R}q(x)dx\to 0$ as $n\to \infty$. Recalling that $q$ is radial, this integral can be written as
\begin{align}\label{integral}
    n\int_{|x|\geq R}q(x)dx=&s_{d-1}na(R)(R)^{d-1}q(R)\\ \nonumber
    &\times \int_0^\infty\left(1+\frac{a(R)}{R}z\right)^{d-1}\frac{q(R+a(R)z)}{q(R)}dz,
\end{align}
where $s_{d-1}$ denotes the surface area of the $(d-1)$ dimensional unit sphere in $\R^d$.\\
Let us estimate the integral on the RHS of (\ref{integral}), as $n\to \infty$. By the mean value theorem, there exists $t$ between $R$ and $R+a(R)z$ such that
$$
\psi(R)-\psi(R+a(R)z)=-(\psi)'(t)(a(R)z).
$$
Since $a(R)=o(R)$, $t\sim R$ and by the regular variation of $a$
$$
(\psi)'(t)=a(t)^{-1}\sim a(R)^{-1};
$$
hence, as $n\to \infty$,
$$
\psi(R)-\psi(R+a(R)z)\to -C'z,
$$
for some constant $C'>0$, from which we find by the dominated convergence theorem that the integral on the RHS above converges as $n\to \infty$ to
$$
\int_0^\infty e^{-C'z}dz<\infty.
$$

The result folllows if the remaining factor on the RHS of (\ref{integral}) tends to $0$, as $n\to \infty$. And indeed we have, as in the proof of Theorem $4.5$ in \cite{topoCrackle}, that
$$
na(R)(R)^{d-1}e^{-\psi(R)}\sim (\log n)^{-1}.
$$
\end{proof}

\section{Subexponential or exponential decay: decrackling the noise}
In this section, suppose that $\psi\in RV_v$ with $v\leq 1$. From now on, let $(r_n)_{n\in \N}$ be a regularly varying sequence decaying to $0$ and such that
\begin{equation}\label{light tail bandwidth}
\frac{a\circ \psi^\leftarrow(\log n)\log\log n}{r_n \psi^{\leftarrow}(\log n)}=o(1).
\end{equation}
Such a choice is feasible by the preliminary observations in Lemma \ref{preliminaries RV}.

\begin{lemma}\label{bandwidth conditions}
It follows from (\ref{light tail bandwidth}) that for $n$ sufficiently large
$$
-\log (r_n)<\log \log n.
$$
\end{lemma}
\begin{proof}
By the above preliminary observations, since here $\psi\in RV_{v}$ with $v\leq 1$, we must have
$$
a\circ\psi^\leftarrow(\log n)\log \log n=\omega(1).
$$
 Using  (\ref{light tail bandwidth}), we then have
\begin{align*}
    (\psi^\leftarrow)^{-1}(\log n) &= o\left(\frac{a\circ\psi^\leftarrow(\log n) \log \log n}{\psi^\leftarrow(\log n)}\right)\\
    &= o(r_n),
\end{align*}
from which it follows, since $r_n=o(1)$, that for $n$ sufficiently large
$$
\log \log n\sim \log\psi^\leftarrow (\log n)>-\log r_n.
$$
\end{proof}
In particular $\log r_n = o(\log n)$ is still true, so the assumptions of Lemma \ref{critical covering} hold.\\

\begin{thm}[Subexponential or exponential decay: decrackling the noise]\label{decrackle conditions}
Suppose that $(r_n)_{n\in \N}$ satisfies (\ref{light tail bandwidth}), then the scaling
$$
x\mapsto \psi^\leftarrow(\log(e+\varphi(x))))
$$
is valid and
$$
\lim_{n\to \infty}\mathbb P\left(\bigcup_{x\in \mathcal P_n}B(x,2r_n\psi^\leftarrow(\log(e+\varphi(x))))\text{ is contractible }\right)=1.
$$
\end{thm}
Note that by regular variation $\psi^\leftarrow(\log(\varphi(R_n^c)))\sim \psi^\leftarrow(\log n)\sim R_n^c$, and since $r_n=o(1)$, we have
$$
r_n\psi^\leftarrow(\log(\varphi(R_n^c)))=o(R_n^c).
$$
Thus the above scaling is already shown to be non-trivial, i.e., to satisfy condition (\ref{condition 2}). We now show that it is also noise-killing, i.e., that it satisfies condition (\ref{condition 1}).\\

\begin{lemma}\label{decrackle lemma}
Assuming (\ref{light tail bandwidth}), we have
$$
\overline{R_n}-R_n^c=o(r_n\psi^\leftarrow\log(\varphi(R_n^c))).
$$
\end{lemma}
\begin{proof}
Using the assumptions on $r_n$, we find by direct computations
\begin{align*}
   \log (\varphi(R_n^c))&\sim \log n + \log r_n - \log \log(r_n^{-1}\psi^\leftarrow(\log n))\\
   &\sim \log n,
\end{align*}
hence by the regular variation properties of $\psi^\leftarrow$,
$$
\psi^\leftarrow(\log(\varphi(R_n^c))))\sim \psi^\leftarrow (\log n).
$$
By the mean value theorem, there exists $t_n\in [A_n,B_n]$ such that
$$
\overline{R_n}-R_n^c=(\psi^\leftarrow)'(t_n)(B_n-A_n),
$$
and ($A_n\sim \log n$ and $B_n\sim \log n)\Rightarrow t_n\sim \log n.$ Using the above preliminary observations, we then find
$$
(\psi^\leftarrow)'(t_n)\sim(\psi^\leftarrow)'(\log n)= a\circ \psi^\leftarrow(\log n).
$$
Thus
$$
(r_n\psi^\leftarrow\log (\varphi(R_n^c)))^{-1}(\overline{R_n}-R_n^c)\sim \frac{a\circ \psi^\leftarrow(\log n)}{r_n}\frac{B_n-A_n}{\psi^\leftarrow(\log n)}.
$$
By Lemma \ref{bandwidth conditions}, we know that for $n$ sufficiently large
$$
-\log r_n < \log \log n,
$$
hence
$$
B_n-A_n\sim \log \log n-\log r_n\sim \log \log n.
$$
Therefore
$$
(r_n\psi^\leftarrow\log (\varphi(R_n^c)))^{-1}(\overline{R_n}-R_n^c)\sim \frac{a\circ\psi^\leftarrow (\log n) \log \log n}{r_n \psi^\leftarrow(\log n)}=o(1).
$$
\end{proof}
We have just shown that the scaling in Theorem \ref{decrackle conditions} is \textit{valid}. It remains to show that the union of the balls is contractible w.h.p..\ Combining Lemmas \ref{critical covering}, \ref{crackle radius} and \ref{decrackle lemma}, this second claim of Theorem \ref{decrackle conditions} holds as follows.

\begin{proof}[Proof of Theorem \ref{decrackle conditions}]
For all $x\in \R^d$ note that
$$
r_n< r_n\psi^\leftarrow(\log (e+\varphi(x)));
$$
hence for $z\in B(0,\overline{R_n}),$ using Lemmas \ref{decrackle lemma} and \ref{critical covering}, and
$$
\dist(z,\mathcal P_n)\leq \dist(z,B(0,R_n^c))+\dist(B(0,R_n^c),\mathcal P_n),
$$
we deduce that with probability going to $1$ as $n\to \infty,$ there exists $x\in \mathcal P_n$ such that
$$
|x-z|<2r_n\psi^\leftarrow(\log(e+\varphi(x))).
$$
Combined with Lemma \ref{crackle radius}, we have with probability going to $1$ as $n\to \infty$
$$
\mathcal P_n\subset B(0,\overline{R_n})\subset \bigcup_{x\in \mathcal P_n}B(x,2r_n\psi^\leftarrow(\log (e+\varphi(x)))),
$$
which gives the theorem.

\end{proof}
\section{Superexponential decay}
In the case where $\psi\in RV_{v}$ and $v>1$, so that the sampling density has superexponential decay, the authors in \cite{topoCrackle} showed that the homology of the noise vanishes provided $r_n$ satisfies the asymptotic condition
$$
\frac{a\circ\psi^\leftarrow(\log n)\log \log n}{r_n}=o(1).
$$
They asked whether this condition could be weakened to merely
$$
\frac{a\circ\psi^\leftarrow(\log n)}{r_n}=o(1).
$$
Here we show that under a suitable \textit{valid} variable bandwidth construction, one may weaken the asymptotic conditions for $r_n$ in this way.

\begin{thm}\label{superexp weakened}
Suppose that
$$
\frac{a\circ\psi^\leftarrow(\log n)}{r_n}=o(1),
$$
and let
\begin{align*}
\mathcal L:\R_+&\longrightarrow \R_+\\
z&\longmapsto \log(e+\log(1+z)).
\end{align*}
The scaling 
$$
x\mapsto \mathcal L (\varphi(x))
$$
is valid and
\begin{equation}\label{superexp contract}
\lim_{n\to \infty}\mathbb P\left(\bigcup_{x\in \mathcal P_n}B(x,2r_n\mathcal L(\varphi(x)))\text{ is contractible }\right)=1.
\end{equation}
\end{thm}
\begin{proof}
We find as before by direct computations,
$$
\log\log \varphi(R_n^c)\sim \log \log n.
$$
Thus we already see that
$$
\mathcal L(\varphi(R_n^c))\sim \log \log n,
$$
while
$$
R_n^c\sim \psi^\leftarrow(\log n),
$$
with $\psi^\leftarrow\in RV_{1/v}$. In particular, we see that the above scaling is non-trivial, i.e., it satisfies condition (\ref{condition 2}).\\

As in the proof of Theorem \ref{decrackle conditions}, by the mean value theorem, there exists\\ $t_n\sim \log n$ such that
$$
\overline{R_n}-R_n^c=(\psi^\leftarrow)'(t_n)(B_n-A_n).
$$
We have as before
$$
(\psi^\leftarrow)'(t_n)\sim(\psi^\leftarrow)'(\log n)= a\circ \psi^\leftarrow(\log n),
$$
hence
$$
(r_n\log\log \varphi(R_n^c))^{-1}(\overline{R_n}-R_n^c)\sim \frac{a\circ \psi^\leftarrow(\log n)}{r_n}\frac{B_n-A_n}{\log \log n}.
$$
The first factor on the RHS above tends to $0$ as $n\to \infty$, by assumption. This also implies that for $n$ sufficiently large
$$
    \log a\circ\psi^\leftarrow(\log n)<\log r_n<0,
$$
hence
$$
-\log r_n \lesssim \log \log n.
$$
Thus, we have
$$
\frac{B_n-A_n}{\log \log n}\lesssim \frac{\log \log n - \log r_n}{\log\log n}=O(1).
$$
Wrapping up, we have just showed that
$$
\overline{R_n}-R_n^c=o(r_n\log \log (\varphi(R_n^c))),
$$
i.e., that the scaling also satisifes condition (\ref{condition 1}), hence is \textit{valid}.\\

The second claim of the theorem, (\ref{superexp contract}),
now follows as in the proof of Theorem \ref{decrackle conditions}.

\end{proof}
\newpage
\section{Conclusion}

This paper is motivated by results in \cite{crackle, topoCrackle}, which investigated conditions under which \v Cech homology of unbounded noise may vanish (which is the case when the union of the balls is contractible, by the Nerve Lemma). In \cite{crackle, topoCrackle} it was shown that unbounded noise introduces non-vanishing homology in general, what the authors call \textit{topological crackle}, unless the sampling density has superexponential decay on $\R^d$.\\

We showed how some well-chosen variable bandwidth constructions allow us to \textit{decrackle} the noise for light tail densities in a non-trivial way. In the case of a density with superexponential decay, which is already known to be exempt from crackle, a well-chosen construction also allows us to weaken the asymptotic conditions on the bandwidth parameter, to address a question posed by the authors in \cite{topoCrackle}.\\

One could naturally be led to ask whether similar results (to the ones obtained above and those contained in \cite{crackle,topoCrackle}) extend to the case where the density is supported on a more general unbounded manifold. In light of the arguments in \cite{crackle,topoCrackle} and those above, such results would only be achievable provided one has a good quantitative estimate for the volume of a Riemannian ball of radius $R(n)\to \infty$, which is highly dependent on the conditions imposed on the manifold. For instance, one could expect some of the results obtained in \cite{topoCrackle} to generalize to the case of an unbounded manifold with non-negative Ricci curvature (via the Bishop-Gromov inequality). It is possible that well-chosen variable bandwidth constructions can help to extend such results to more general unbounded manifolds.

\newpage

\bibliographystyle{plain}
\bibliography{refs}
\end{document}